\documentclass[10pt,twoside]{article}
\usepackage[english]{babel}
\usepackage{droid}
\usepackage{famems}
\usepackage{braket}
\usepackage{citeref}
\usepackage{cmap}
\usepackage{graphics}
\usepackage{graphicx}
\usepackage{graphpap}
\usepackage{amsmath,amstext,amsfonts,amssymb,amsthm}
\usepackage{mathtools}
\usepackage{dsfont}
\usepackage[mathscr]{eucal}
\usepackage{array}
\usepackage{multicol}
\usepackage{multirow}
\usepackage{multicol}
\usepackage{tabularx}
\usepackage{xcolor}
\usepackage{color}
\usepackage{colortbl}
\usepackage{ifthen}
\usepackage{pdflscape}
\usepackage{epstopdf}
\usepackage[normalem]{ulem}
\usepackage{rotating}
\usepackage{hhline}
\usepackage{fontawesome}
\usepackage{scrextend}
\usepackage{arydshln}
\usepackage{ctable}
\usepackage{url}
\usepackage[breaklinks]{hyperref}
\usepackage{hyperref}
\hypersetup{
    colorlinks,
    citecolor=black,
    filecolor=black,
    linkcolor=black,
    urlcolor=
teal
}
\usepackage{cleveref}
\definecolor{MyYellow50}{RGB}{255,255,128}
\definecolor{MyYellow0}{RGB}{255,255,0}
\definecolor{MyYellow75}{RGB}{255,255,192}
\definecolor{MyYellow85}{RGB}{255,255,240}
\definecolor{MyAqua50}{RGB}{127,255,255}
\definecolor{MyAqua0}{RGB}{0,255,255}
\definecolor{MyAqua75}{RGB}{192,255,255}
\definecolor{MyAqua85}{RGB}{224,255,255}
\definecolor{MyGreen0}{RGB}{0,255,0}
\definecolor{MyGreen50}{RGB}{128,255,128}
\definecolor{MyYellow}{rgb}{1,1,0}
\definecolor{MyYellow75}{rgb}{1,1,0.75}
\definecolor{MyRed95}{rgb}{1,0.95,0.95}
\definecolor{MyRed85}{rgb}{1,0.85,0.85}
\definecolor{MyRed75}{rgb}{1,0.75,0.75}
\definecolor{MyRed50}{rgb}{1,0.50,0.50}
\definecolor{MyRed25}{rgb}{1,0.25,0.25}
\definecolor{MyRed0}{rgb}{1,0,0}
\definecolor{MyRed}{rgb}{1,0,0}
\definecolor{MyBlue}{rgb}{0,0,1}
\definecolor{MyGreen}{rgb}{0.75,1,0.75}
\definecolor{MyGreen35}{rgb}{0.35,1,0.35}
\definecolor{MyGreen75}{rgb}{0.75,1,0.75}
\definecolor{MyBlue40}{rgb}{0.40,0.40,1}
\definecolor{MyBlue60}{rgb}{0.60,0.60,1}
\definecolor{MyBlue75}{rgb}{0.75,0.75,1}
\definecolor{MyBlue85}{rgb}{0.85,0.85,1}
\definecolor{MyGray}{rgb}{0.95,0.95,0.95}
\definecolor{MyWhite}{rgb}{1,1,1}
\definecolor{MyBlack}{rgb}{0,0,0}
\definecolor{MyProbability}{rgb}{1.00,0.25,0.25}
\definecolor{MyConfidency}{rgb}{1.00,0.25,1.00}
\definecolor{MyEventology}{rgb}{0.25,1.00,1.00}
\definecolor{MyOrange}{rgb}{1,0.56,0.25}
\definecolor{MyOrangeB}{rgb}{1,0.65,0.45}
\definecolor{MyMagentaP}{rgb}{1.00,0.45,0.65}
\definecolor{MyMagenta}{rgb}{1.00,0,1.00}
\definecolor{MyAqua}{rgb}{0,0.75,0.75}
\makeatletter
\renewcommand{\paragraph}{\@startsection{paragraph}{4}{0ex}%
   {-3.25ex plus -1ex minus -0.2ex}%
   {1.5ex plus 0.2ex}%
   {\normalfont\normalsize\tt}}
\makeatother

\begin{document}
\baselineskip=11pt

\newcounter{ctrwar}\setcounter{ctrwar}{0} 
\newcounter{ctrdef}\setcounter{ctrdef}{0}
\newcounter{ctrdefpre}\setcounter{ctrdefpre}{0}
\newcounter{ctrTh}\setcounter{ctrTh}{0}
\newcounter{ctrnot}\setcounter{ctrnot}{0}
\newcounter{ctrATT}\setcounter{ctrATT}{0} 
\newcounter{ctrcor}\setcounter{ctrcor}{0}
\newcounter{ctrAx}\setcounter{ctrAx}{0}
\newcounter{ctrexa}\setcounter{ctrexa}{0}
\newcounter{ctrlem}\setcounter{ctrlem}{0}
\newcounter{ctrPRO}\setcounter{ctrPRO}{0}
\newcounter{ctrrem}\setcounter{ctrrem}{0}
\newcounter{ctrass}\setcounter{ctrass}{0}
\newcounter{ctrmem}\setcounter{ctrmem}{0}

\newcommand{\bfPB}{\mbox{\protect\reflectbox{\bf P}\hspace{-0.4em}{\bf B}}}
\newcommand{\bfEl}{\mbox{\protect\reflectbox{\bf E}}}
\newcommand{\bfEr}{\mathbf{E}}
\newcommand{\bfEE}{\mbox{\protect\reflectbox{\bf E}\hspace{-0.4em}{\bf E}}}
\newcommand{\bfphi}{\mbox{\boldmath$\varphi$}}
\newcommand{\bfPhi}{\mbox{\boldmath$\Phi$}}
\newcommand{\bfPhii}{\mbox{\scriptsize\boldmath$\Phi$}}
\newcommand{\rS}{\reflectbox{\bf S}}
\newcommand{\rsS}{\,\reflectbox{\scriptsize\bf S}}
\newcommand{\rssS}{\,\reflectbox{\tiny\bf S}}
\newcommand{\reS}{\reflectbox{S}}
\newcommand{\resS}{\,\reflectbox{\scriptsize S}}
\newcommand{\ressS}{\,\reflectbox{\tiny S}}
\newcommand{\scrER}{\mathscr E\!\mathscr R}

\newcommand{\ER}{\mathscr E\!\mathscr R}

\numberwithin{equation}{section}

\renewcommand{\figurename}{\scriptsize  Figure}
\renewcommand{\tablename}{\scriptsize  Table}
\renewcommand\refname{\scriptsize References}
\renewcommand\contentsname{\scriptsize Contents}
\setcounter{page}{63}

\titleOneAuthoren{
Fr\'echet bounds of the 1-st kind\\[-3pt] for sets of half-rare events

}{Oleg Yu. Vorobyev}{
Institute of Mathematics and Computer Science\\
Siberian Federal University\\
Krasnoyarsk\\
\tiny
\url{mailto:oleg.yu.vorobyev@gmail.com}\\
\url{http://www.sfu-kras.academia.edu/OlegVorobyev}
}

\label{vorobyev-17}

\colontitleen{Vorobyev}
\footavten{Oleg Yu. Vorobyev}
\setcounter{footnote}{0}
\setcounter{equation}{0}
\setcounter{figure}{0}
\setcounter{table}{0}
\setcounter{section}{0}

\vspace{-22pt}

\begin{abstracten}
\textit{Fr\'echet bounds of the 1-st kind for sets of events (s.e.) and its main properties are considered. The lemma on not more than two nonzero values of lower Fr\'echet-bounds of the 1-st kind for a set of half-rare events (s.h-r.e.) is proved with the corollary on the analogous assertion for s.e's with arbitrary event-probability distributions (e.p.d's.).}
\end{abstracten}\\[-21pt]
\begin{keywordsen}
\emph{Probability, Kolmogorov event, event, set of events, event-probability distribution, set of half-rare events, Frechet bounds of the 1st kind.}
\end{keywordsen}


In \cite{Vorobyev2015famems12} it is shown that any s.e's are set-phenomena of some \emph{half-rare s.e.} And its e.p.d's are characterized by the e.p.d's of the half-rare s.e's the set-phenomena of which they are. In view of this, a theory of Fr\'echet-bounds of the 1st kind of \emph{half-rare s.e's} is being developed in this paper, i.e., the theory of such s.e's events of which happen with probabilities not greater than half. These ``nondescript'' restrictions are enough to build Fr\'echet-bounds for any s.e. from Fr\'echet-bounds of the half-rare s.e., relying on set-phenomenon transformations of e.p.d's and of its Fr\'echet-bounds of the 1st kind.

We will make one more purely technical assumption, which does not detract from the generality of our analysis. We will assume that events in each half-rare s.e. are ordered in descending of their probabilities and we introduce the following natural notation and abbreviations for half-rare $N$-s.e.
$
\frak{X} = \{ x, x \in \frak{X} \} = \{ x_1,...,x_N \}
$
and its $\frak{X}$-set of probabilities of events
$
\breve{p} = \{p_x, x \in \frak{X}\} = \{p_{x_1},...,p_{x_N}\} = \{p_1,...,p_N\}
$
where
$$
1/2 \geqslant p_1 \geqslant p_2 \geqslant ... \geqslant p_N.
\eqno{(\star)}
$$

\section{Properties of Fr\'echet-bounds of the 1st kind of a half-rare s.e.\label{FBproperties}}

Consider the half-rare $N$-s.e. $\frak{X}=\{x_1,...,x_N\}$ with the $\frak{X}$-set of probabilities of events $\breve{p}=\{p_x, x \in \frak{X}\}=\{p_1,...,p_N\}$ where, we recall,
$
1/2 \geqslant p_1 \geqslant ... \geqslant p_N,
$
and the e.p.d. of the 1st kind
$
\{p(X /\!\!/ \frak{X}), X \subseteq \frak{X}\},
$
and Fr\'echet-bounds of the 1st kind for $X \subseteq \frak{X}$:
$$
p^-(X /\!\!/ \frak{X}) \leqslant p(X /\!\!/ \frak{X}) \leqslant p^+(X /\!\!/ \frak{X})
$$
where
$$
p^-(X /\!\!/ \frak{X})
= \max\left\{ 0, 1-\hspace{-5pt}\sum_{x \in X} \! (1-p_x)-\hspace{-10pt}\sum_{x \in \frak{X}-X} \!\!\!\! p_x \right\}\!\!,
\eqno{(\ref{FBproperties}.1)}
$$
$$
p^+(X /\!\!/ \frak{X}) = \min \left\{ \min_{x \in X} p_x, \, \min_{x \in \frak{X}-X} (1-p_x) \right\}
\eqno{(\ref{FBproperties}.2)}
$$
are general formulas for the lower and the upper Fr\'echet-bounds of the 1st kind. They form the sets of Fr\'echet-bounds which are called the lower and the upper \emph{Fr\'echet-boundary distributions (F.-b.d's)} of the 1st kind:
$$
\{p^-(X /\!\!/ \frak{X}), X \subseteq \frak{X}\}, \ \ \ \{p^+(X /\!\!/ \frak{X}), X \subseteq \frak{X}\}.
$$

The upper Fr\'echet-bound can be simplified by virtue of the characteristic property of a half-rare s.e. ($\star$) since then for non-empty $\emptyset \not= X \subseteq \frak{X}$
$$
\min_{z \in X} p_z \leqslant \min_{z \in \frak{X}-X} (1-p_z).
$$
We have the formula for non-empty $\emptyset \not= X \subseteq \frak{X}$:
$$
p^+(X /\!\!/ \frak{X}) = \min_{z \in X} p_z.
$$
If $X=\emptyset$, then, as is not difficult to understand,
$$
p^+(\emptyset /\!\!/ \frak{X}) = \min_{z \in \frak{X}} (1-p_z) = 1-p_x.
$$
From this the formula for the upper Fr\'echet-bound of the 1st kind of a half-rare s.e. has the form:
$$
p^+(X /\!\!/ \frak{X}) =
\begin{cases}
1-p_x,              & X=\emptyset,\cr
\displaystyle \min_{z \in X} p_z, & \emptyset \not= X \subseteq \frak{X}.
\end{cases}
\eqno{(\ref{FBproperties}.3)}
$$
and the \emph{upper F.-b.d. of the 1st kind} has the form:
$$
\left\{1-p_x, \min_{z \in X} p_z, \emptyset \not= X \subseteq \frak{X}\right\}.
$$

With the lower Fr\'echet-bounds of the 1st kind of any half-rare s.e. $\frak{X}$ the situation is not much more complicated.

\textbf{Lemma \refstepcounter{ctrlem}\arabic{ctrlem}\label{lem1}\itshape\enskip\scriptsize (on the lower Fr\'echet-bounds of the 1st kind of a half-rare s.e.).} The \emph{lower Fr\'echet-bound of the 1st kind of the half-rare s.e. $\frak{X}$ with maximum probability of events
$\displaystyle p_{\max}=\max_{z \in \frak{X}} p_z = \mathbf{P}(x_{\max})$ can be nonzero for only two probabilities of the 1st kind: $p(\emptyset /\!\!/ \frak{X})$ and $p(\{x_{\max}\} /\!\!/ \frak{X})$}.

\texttt{Proof.}
From the definition of the lower Fr\'echet-bound of the 1st kind of the s.e. $\frak{X}$ (\ref{FBproperties}.1)
and due to the fact that $\frak{X}$ is the half-rare s.e., it is clear that for $X \subseteq \frak{X}$ such that $|X|>1$, it is equal to zero, since in this case there are two such events $x, x' \in X$ that $1-p_x>1/2$ and $1-p_{x'}>1/2$. An therefore
$\displaystyle
1-\sum_{x \in X} (1-p_x) < 0
$
and $p^-(X /\!\!/ \frak{X}) =0$.
Consequently, for the arbitrary half-rare s.e. $\frak{X}$ the lower Fr\'echet-bound of the 1st kind can differ from zero only for such $X \subseteq \frak{X}$ the power of which is not greater than unity --- $|X| \leqslant 1$:
$$
p^-(X /\!\!/ \frak{X}) =
\begin{cases}
p^-(X /\!\!/ \frak{X}) \geqslant 0, & |X| \leqslant 1, X \subseteq \frak{X},\cr
0, & |X| > 1, X \subseteq \frak{X}.
\end{cases}
$$
Since for the half-rare $N$-s.e. $\frak{X}$ there is only polynomial number, $N+1$, of such subsets then задача the problem of constructing the lower Fr\'echet-bound of the 1st kind of a half-rare s.e. we can consider basically solved. After detailing the condition $|X|\leqslant 1$ the formula takes the following form:
$$
p^-(X /\!\!/ \frak{X}) =
$$
$$
=
\begin{cases}
\displaystyle \max\left\{0, \, 1-\sum_{x \in \frak{X}} p_x \right\}, & X = \emptyset,\cr
\displaystyle \max\left\{0, \, p_x-\sum_{z \in \frak{X}-\{x\}} p_z \right\}, & x \in \frak{X},\cr
0, & |X| > 1, X \subseteq \frak{X}.
\end{cases}
$$
Now is just the time, once again take advantage of the fact that $\frak{X}$ is a half-rare s.e. Quite unexpectedly it turns out that for the half-rare s.e. $\frak{X}$ among $N$ conditions for $x \in \frak{X}$ only for the event $x_{\max} \in \frak{X}$,
that happens with maximum probability
$\displaystyle p_{\max} = \max_{z \in \frak{X}} p_z$,
the expression
$$
p_{\max}-\sum_{z \in \frak{X}-\{x_{\max}\}} p_z
$$
under the sign of maximum can be strictly greater than zero. For the rest events $z \not= x_{\max}$ from $\frak{X}$ this expression can be only not greater than zero, since from its non-maximum probability $p_z$,
since the maximum probability  and all the rest ones will be subtracted from it:
$$
p_z-p_{\max}-\sum_{y \in \frak{X}-\{x_{\max},z\}} p_y \leqslant 0.
$$
This implies the assertion of the lemma, since we obtain the following formula for the lower Fr\'echet-bound of the 1st kind of a half-rare s.e. for $X \subseteq \frak{X}$:
$$
p^-(X /\!\!/ \frak{X}) =
$$
$$
=
\begin{cases}
\displaystyle \max\left\{0, \, 1-\sum_{x \in \frak{X}} p_x \right\}, &\hspace{-8pt} X = \emptyset,\cr
\displaystyle \max\left\{0, \, p_{\max}-\hspace{-10pt}\sum_{z \in \frak{X}-\{x_{\max}\}}\hspace{-10pt} p_z \right\}\!, &\hspace{-8pt}  X=\{x_{\max}\},\cr
0, &\hspace{-8pt} \mbox{иначе},
\end{cases}
\eqno{(\ref{FBproperties}.4)}
$$
and the \emph{lower F.-b.d. of the 1st kind} has the form:
$$
\left\{ p^-(\emptyset /\!\!/ \frak{X}),p^-(\{x_{\max}\} /\!\!/ \frak{X}),0,...,0 \right\}
$$
where two possibly ono-zero the lower Fr\'echet-bounds are defined by the formula (\ref{FBproperties}.4).
The lemma is proved.

\textbf{Corollary \refstepcounter{ctrcor}\arabic{ctrcor}\label{cor1}\itshape\enskip\scriptsize (on the lower Fr\'echet-bounds of the 1st kind of an arbitrary s.e.).} \emph{The assertion of the lemma  is valid not only for half-rare, but also for arbitrary s.e's.}

\texttt{Proof} follows from the lemma on the characterization of an arbitrary s.e. by its half-rare projection \cite{Vorobyev2015famems12}, according to which the e.p.d. of the 1st kind of an arbitrary s.e., as well as its F.-b.d. of the 1st kind are \emph{set-phenomenal renumbering} \cite{Vorobyev2015famems12} of corresponding e.p.d's and F.-b.d's of its half-rare projection.

\textbf{Corollary \refstepcounter{ctrcor}\arabic{ctrcor}\label{cor2}\itshape\enskip\scriptsize (on the lower Fr\'echet-bounds of the 1st kind of a half-rare doublet of events).} \emph{Fr\'echet-bounds of the 1st kind of the half-rare doublet of events $\frak{X}=\{x,y\}$ with the $\{x,y\}$-set of marginal probabilities $\breve{p}=\{p_x,p_y\}$, such that $1/2 \geqslant p_x \geqslant p_y$, have the form:}
\begin{equation}\nonumber
\hspace{10pt}\begin{split}
1-p_x-p_y = p^-(\emptyset) \leqslant p(\emptyset) &\leqslant p^+(\emptyset) = 1-p_x,\\
p_x-p_y = p^-(x) \leqslant p(x) &\leqslant p^+(x) = p_x,\\
0 = p^-(y) \leqslant p(y) &\leqslant p^+(y) = p_y,\\
0 = p^-(xy) \leqslant p(xy) &\leqslant p^+(xy) = p_y.
\end{split}
\hspace{0pt}(\ref{FBproperties}.5)
\end{equation}
\texttt{Proof}. The formulas (\ref{FBproperties}.5) are an immediate consequence of formulas (\ref{FBproperties}.3) and (\ref{FBproperties}.4).

\textbf{Corollary \refstepcounter{ctrcor}\arabic{ctrcor}\label{cor3}\itshape\enskip\scriptsize (on Fr\'echet-inequalities for a covariance of the 1st kind of a half-rare doublet of events).} \emph{Fr\'echet-inequalities for a covariance of the 1st kind of the half-rare doublet of events $\frak{X}=\{x,y\}$ with the $\{x,y\}$-set of marginal probabilities $\breve{p}=\{p_x,p_y\}$, such that $1/2 \geqslant p_x \geqslant p_y$, have the form:}
\begin{equation}\nonumber
\hspace{15pt}\begin{split}
-p_x p_y &\leqslant \textsf{Kov}(\emptyset) \leqslant (1-p_x)p_y,\\
-(1-p_x)p_y &\leqslant \textsf{Kov}(x) \leqslant p_x p_y,\\
-(1-p_x)p_y &\leqslant \textsf{Kov}(y) \leqslant p_x p_y,\\
-p_x p_y &\leqslant \textsf{Kov}(xy) \leqslant (1-p_x)p_y.
\end{split}
\hspace{30pt}(\ref{FBproperties}.6)
\end{equation}
\texttt{Proof} follows from the definition of covariance of the 1st kind and Corollary \ref{cor2}.

\section{The pair illustrations\label{FBillustrations}}

It is interesting to look at the graphs of F.-b.d's of the 1st kind, in order to ``make sure with our own eyes'' in the just proven of their amazing properties.

\begin{figure}[h!]
\vspace*{-10pt}
\centering
\includegraphics[width=1.6in]{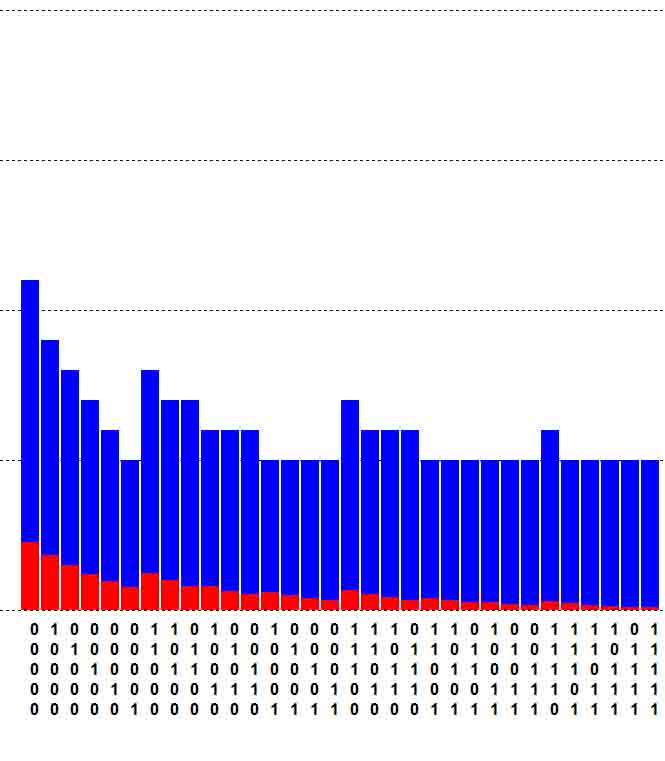}  

\vspace{-5pt}

\includegraphics[width=0.125in]{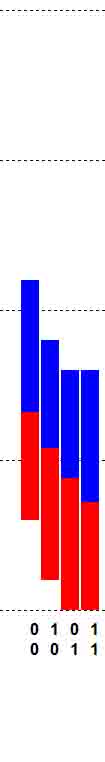} 
\includegraphics[width=0.213in]{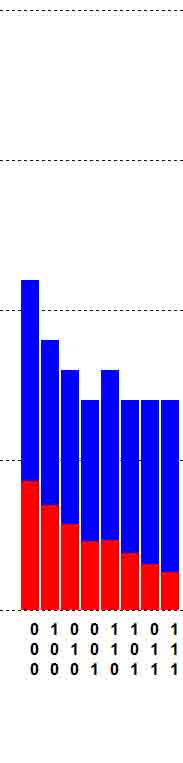} 
\includegraphics[width=0.40in]{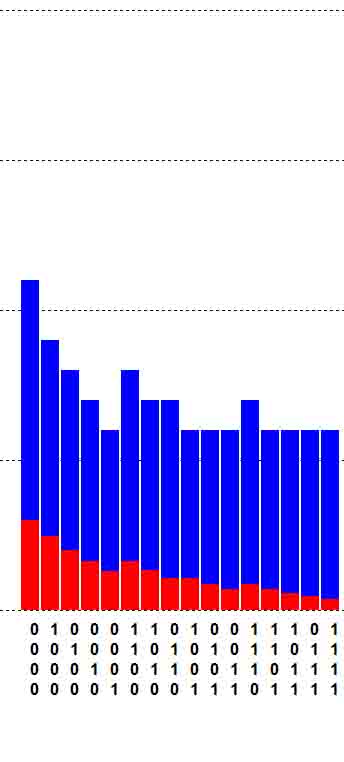} 
\includegraphics[width=0.775in]{pentapletDOWN.jpg}  
\includegraphics[width=1.53in]{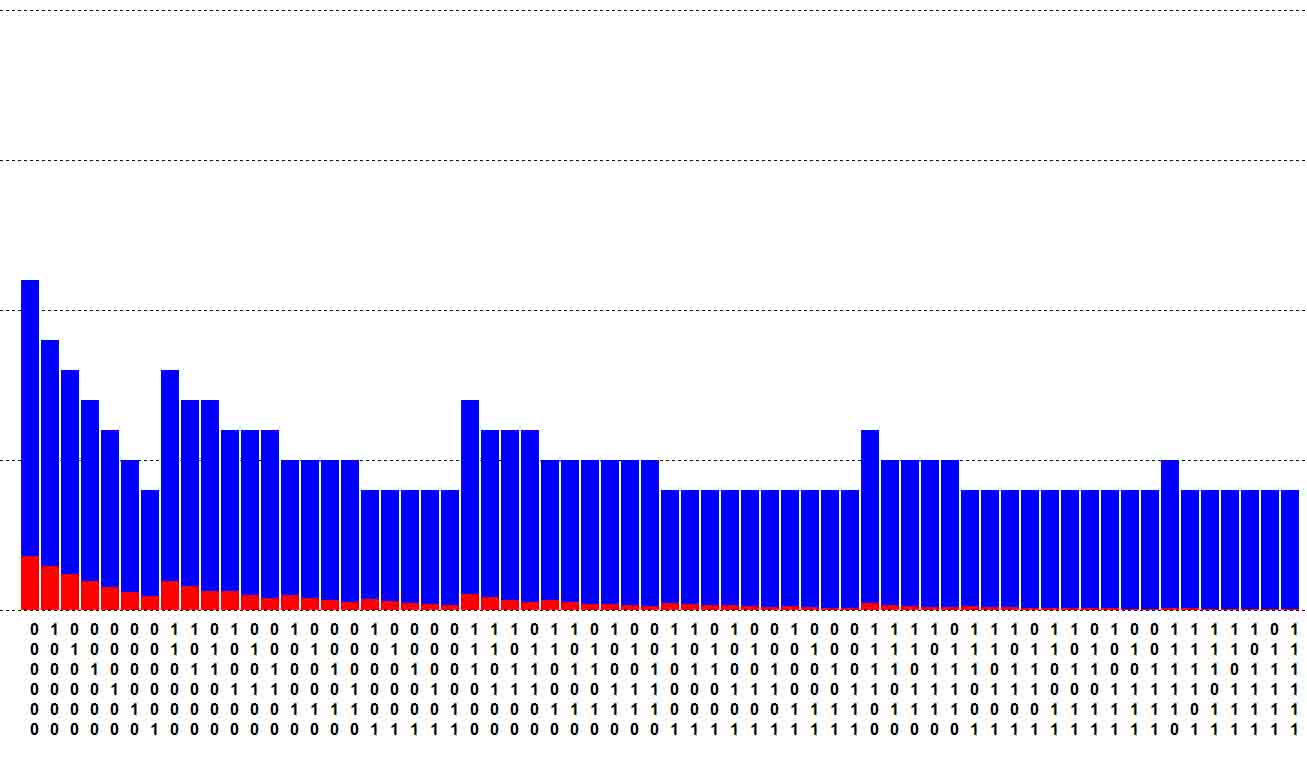} 

\includegraphics[width=3.34in]{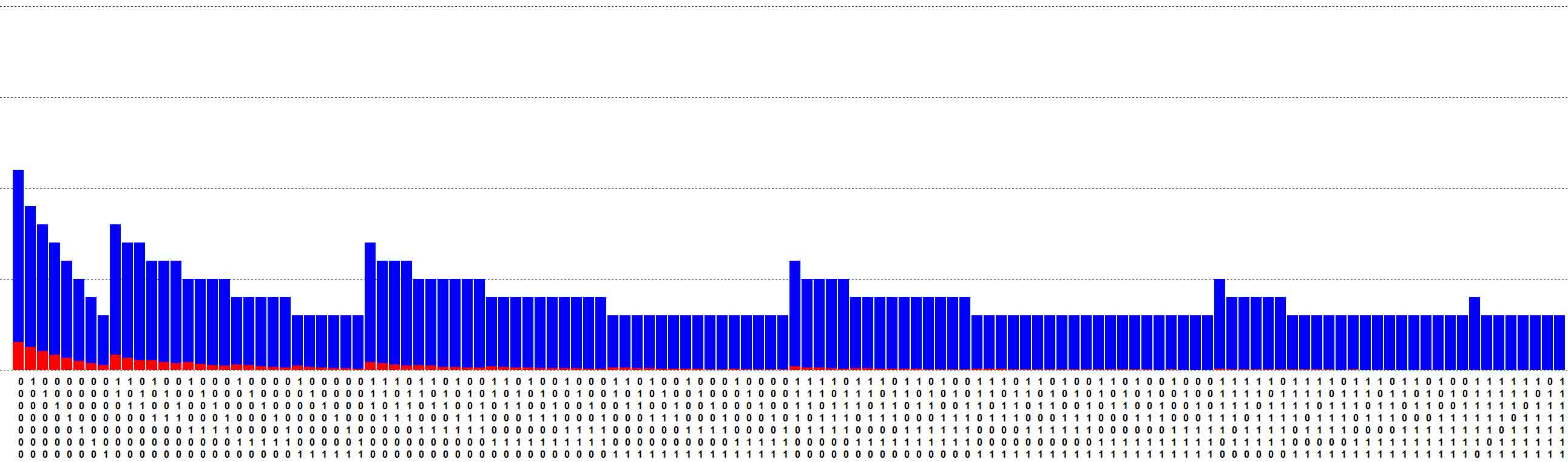} 

\vspace{-5pt}

\caption{Standard s.e's $\frak{X}$ of small power: $N=2,3,4,5,6,7$ with $\frak{X}$-sets of probabilities of events from $\breve{p}=\{0.45,0.40\}$ for a doublet to $\breve{p}=\{0.45,0.40,0.35,0.30,0.25,0.20,0.15\}$ for septa-plet;  At the top of the close-up is shown again the graph for the half-rare penta-plet  $\frak{X}$ with the $\frak{X}$-set of probabilities of event $\breve{p}=\{0.45,0.40,0.35,0.30,0.25\}$. \label{figsDOWN}} %
\end{figure}

One of the concepts that appear on these graphs is the concept of the \emph{independent projection of an arbitrary s.e.} $\frak{X}$, by which, recall, we understood a s.e., e.p.d. of the 1st kind of which is called \emph{independent e.p.d. (i.e.p.d.) of the 1st kind}; it is denoted by
$$
\{p^\star(X /\!\!/ \frak{X}), X \subseteq \frak{X})\}
$$
and it is defined for $X \subseteq \frak{X}$ by formulas
$$
p^\star(X /\!\!/ \frak{X}) = \prod_{x \in X} p_x \prod_{x \in \frak{X}-X} (1-p_x)
$$
where $p_x = \mathbf{P}(x), x \in \frak{X}$, are probabilities of events from $\frak{X}$.

\begin{figure}[h!]
\vspace*{-10pt}
\centering
\includegraphics[width=1.05in]{pentapletDOWN.jpg}  
\includegraphics[width=1.05in]{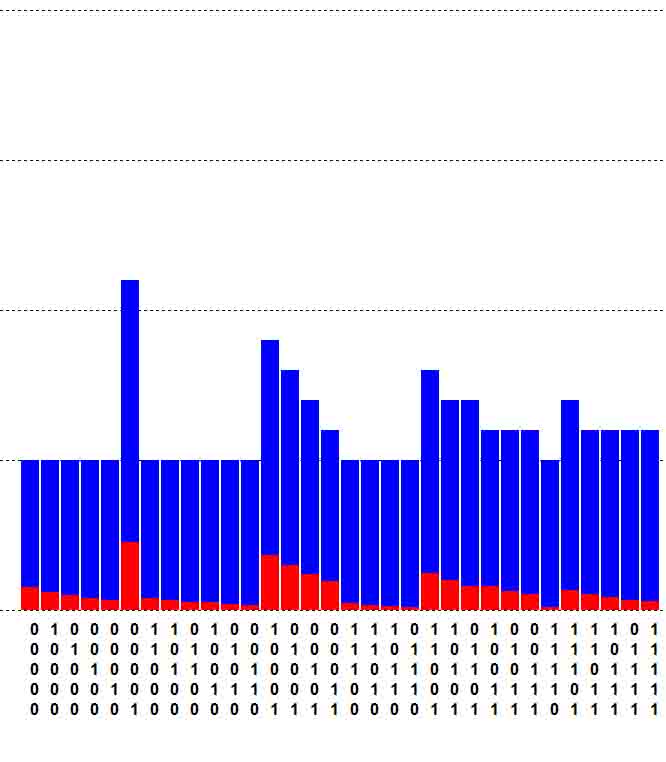}  
\includegraphics[width=1.05in]{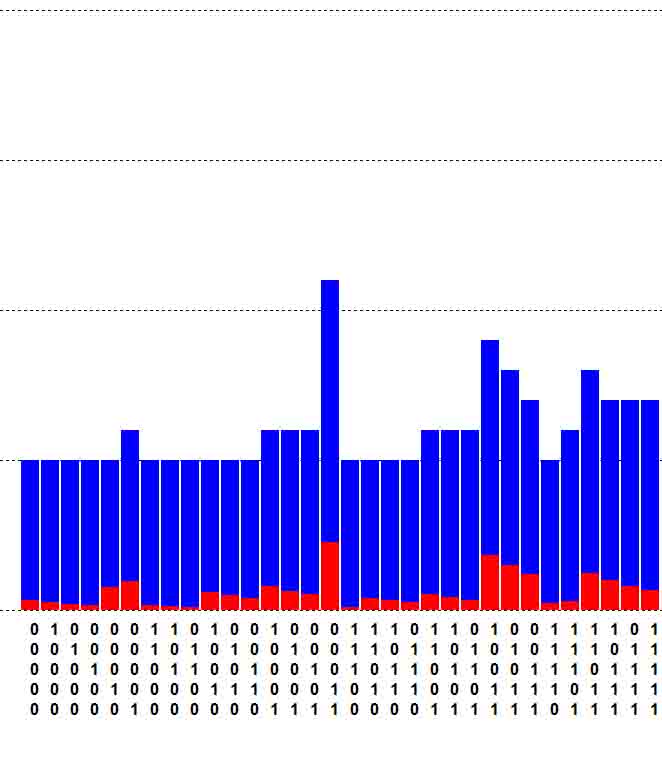}  

\includegraphics[width=1.05in]{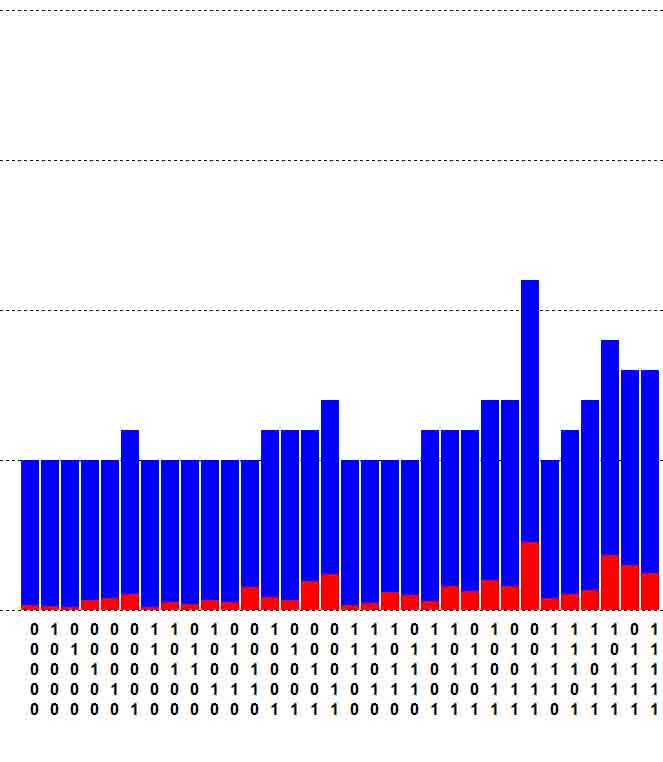}  
\includegraphics[width=1.05in]{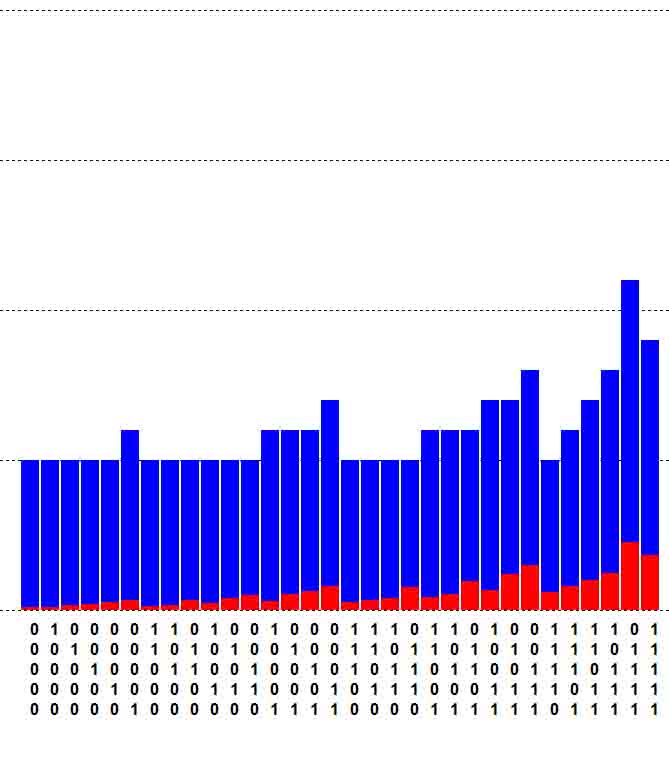}  
\includegraphics[width=1.05in]{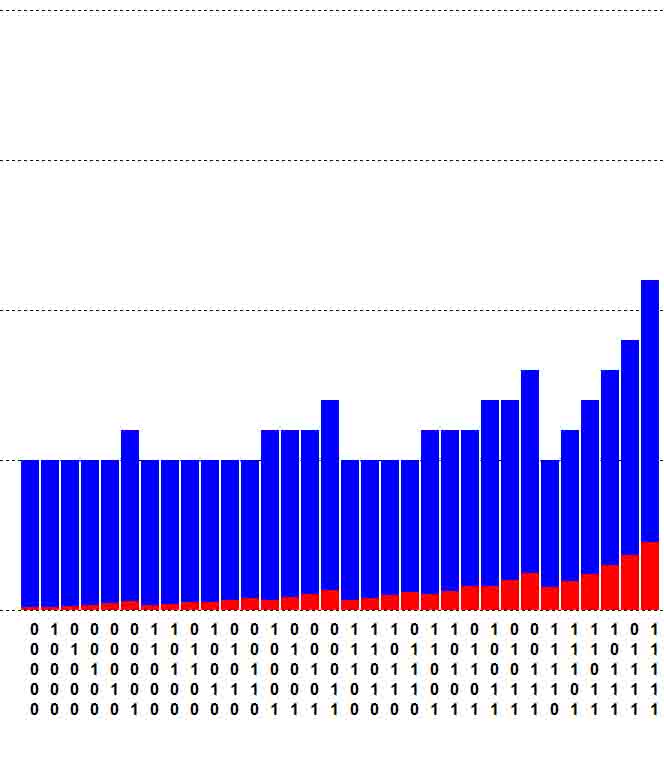}  

\vspace{-5pt}

\caption{\textbf{Standard penta-plet} with the $\frak{X}$-sets of probabilities of events $\breve{p}=\{0.45,0.40,0.35,0.30,0.25\}$ and \textbf{five its set-phenomena}: $\{x_1,x_2,x_3,x_4\}$-phenomenon, $\{x_1,x_2,x_3\}$-phenomenon, $\{x_1,x_2\}$-phenomenon, $\{x_1\}$-phenomenon and $\emptyset$-phenomenon (M-complement).
\label{figsDOWNph}} %
\end{figure}

Each of the graphs illustrates the lower and the upper F.-b.d's and i.e.p.d's of independent projection of the half-rare s.e. $\frak{X}$ of low power: $N=2,3,4,5,6,7$ for $X \subseteq \frak{X}$ and different values of the $\frak{X}$-set of probabilities of events $\breve{p}=\{p_x, x \in \frak{X}\}$. For $X \subseteq \frak{X}$ the interval
\textcolor{blue}{$$
[\,p^\star(X /\!\!/ \frak{X}), \, p^+(X /\!\!/ \frak{X})\,]
$$}
between values of the i.e.p.d. and the upper F.-b.d. is shown as \textcolor{blue}{blue}, and the interval
\textcolor{red}{$$
[\,p^-(X /\!\!/ \frak{X}), \,p^\star(X /\!\!/ \frak{X})\,]
$$}between values of the lower F.-b.d. and the i.e.p.d. is shown as \textcolor{red}{red}.

Subsets of events $X \subseteq \frak{X}$ is denoted by $\frak{X}$-sets:
$
X \sim \{ \mathbf{1}_X(x), x \in \frak{X} \}
$
of values of its indicators on $\frak{X}$. For example, the empty subset $\emptyset \subseteq \frak{X}$ is denoted by the $\frak{X}$-set from zeros:
$
\emptyset \sim \{0,...,0\},
$
and the s.e. $\frak{X}$ by the $\frak{X}$-set from units:
$
\frak{X} \sim \{1,...,1\}.
$
The horizontal dotted line indicates the scale of the unit interval $[0,1]$ along the vertical axis in 1/4.

$\bigstar$ The English version of this article was published on October 2, 2017. Therefore, my later works \cite{Vorobyev2015famems13,Vorobyev2016famems3,Vorobyev2016famems2,Vorobyev2016famems1}, which expand the themes of this work, are added to the list of references.

{\footnotesize
\bibliography{vorobyev5}
}

\end{document}